\pdfoutput=1
\RequirePackage{ifpdf}
\ifpdf 
\documentclass[pdftex]{sigma}
\else
\documentclass{sigma}
\fi

\input xy			
\xyoption{all}
\xyoption{2cell}
\UseComputerModernTips
\UseTwocells

\numberwithin{equation}{section}
\newtheorem{Theorem}{Theorem}[section]
\newtheorem{Corollary}[Theorem]{Corollary}
\newtheorem{Lemma}[Theorem]{Lemma}
\newtheorem{Proposition}[Theorem]{Proposition}
{\theoremstyle{definition}
\newtheorem{Definition}[Theorem]{Definition}
\newtheorem{Example}[Theorem]{Example}
\newtheorem{Remark}[Theorem]{Remark}
}

\DeclareMathOperator{\Spinc}{Spin^{c}}
\DeclareMathOperator{\rank}{rank}

\DeclareMathOperator{\SO}{SO}

\DeclareMathOperator{\SU}{SU}
\DeclareMathOperator{\PU}{PU}

\begin{document}

\newcommand{\refl}[1]{s_{#1}}

\newcommand{\arXivNumber}{1408.0253}

\allowdisplaybreaks

\renewcommand{\thefootnote}{$\star$}

\renewcommand{\PaperNumber}{109}

\FirstPageHeading

\ShortArticleName{Prequantization of the Moduli Space of Flat $\PU(p)$ Bundles}

\ArticleName{Prequantization of the Moduli Space\\
of Flat $\boldsymbol{\PU(p)}$-Bundles\\
with Prescribed Boundary Holonomies\footnote{This paper is a~contribution to the Special Issue on Poisson Geometry in
Mathematics and Physics.
The full collection is available at \href{http://www.emis.de/journals/SIGMA/Poisson2014.html}
{http://www.emis.de/journals/SIGMA/Poisson2014.html}}}

\Author{Derek KREPSKI}

\AuthorNameForHeading{D.~Krepski}

\Address{Department of Mathematics, University of Manitoba, Canada}
\Email{\href{mailto:Derek.Krepski@umanitoba.ca}{Derek.Krepski@umanitoba.ca}}
\URLaddress{\url{http://server.math.umanitoba.ca/~dkrepski/}}

\ArticleDates{Received August 05, 2014, in f\/inal form November 28, 2014; Published online December 05, 2014}

\Abstract{Using the framework of quasi-Hamiltonian actions, we compute the obstruction to prequantization for the moduli
space of f\/lat ${\rm PU}(p)$-bundles over a~compact orientable surface with prescribed holonomies around boundary components,
where $p>2$ is prime.}

\Keywords{quantization; moduli space of f\/lat connections; parabolic bundles}

\Classification{53D50; 53D30}

\renewcommand{\thefootnote}{\arabic{footnote}}
\setcounter{footnote}{0}

\section{Introduction}

Let~$G$ be a~compact connected simple Lie group and~$\Sigma$ a~compact oriented surface with~$s$ boundary components.
Given conjugacy classes $\mathcal{C}_1, \ldots, \mathcal{C}_s$, let $\mathfrak{M}=M_G(\Sigma; \mathcal{C}_1, \ldots,
\mathcal{C}_s) $ denote the moduli space of f\/lat $G$-bundles on~$\Sigma$ with prescribed boundary holonomies in the
conjugacy classes $\mathcal{C}_j$.
Alternatively, $\mathfrak{M}$ may be described as the character variety of the fundamental group of~$\Sigma$,
\begin{gather*}
\mathfrak{M} = \mathrm{Hom}_{\mathcal{C}_1, \ldots, \mathcal{C}_s} (\pi_1({\Sigma}),G)/G.
\end{gather*}
Here, $\mathrm{Hom}_{\mathcal{C}_1, \ldots, \mathcal{C}_s} (\pi_1({\Sigma}),G)$ consists of homomorphisms
$\rho:\pi_1({\Sigma}) \to G$ whose restriction to (the homotopy class of) the~$j$-th boundary circle of~$\Sigma$ lies in
$\mathcal{C}_j$, and~$G$ acts by conjugation.
Recall that $\mathfrak{M}$ is a~(possibly singular) symplectic space, where the symplectic form is def\/ined by a~choice
of invariant inner product on the Lie algebra $\mathfrak{g}$ of~$G$~\cite{AB,Go}.
This paper considers the obstruction to the existence of a~prequantization of $\mathfrak{M}$~-- that is, a~prequantum
(orbifold) line bundle $L\to \mathfrak{M}$ (see Section~\ref{pqbundle} for details)~-- by expressing the corresponding
integrality condition on the symplectic form in terms of the choice of inner product on the simple Lie algebra
$\mathfrak{g}$, which is hence a~certain multiple~$k$ of the basic inner product.

If the underlying structure group~$G$ is simply connected, the moduli space $\mathfrak{M}$ is connected and the
obstruction to prequantization is well known~-- a~prequantization exists if and only if $k\in {\mathbb{N}}$ and each
conjugacy class $\mathcal{C}_j$ corresponds to a~level~$k$ weight (e.g., see~\cite{AMW, BL, BasicGerbe}).
If~$G$ is not simply connected, $\mathfrak{M}$ may have multiple components.
Moreover integrality of~$k$ is not suf\/f\/icient to guarantee a~prequantization even in the absence of markings/prescribed
boundary holonomies: if~$\Sigma$ is closed and has genus at least~$1$, then~$k$ must be a~multiple of an integer~$l_0(G)$ (computed in~\cite{JGP} for each~$G$).
If~$\Sigma$ has boundary with prescribed holonomies, only the case $G=\SO(3)\cong \PU(2)$ has been fully
resolved~\cite{KM}.

In this paper, we describe the connected components of $\mathfrak{M}$ for non-simply connected structure groups $G/Z$ in
Corollary~\ref{cor:comps} and Proposition~\ref{prop:comps} (where~$G$ is simply connected and~$Z$ is a~subgroup of the centre of~$G$).
The decomposition into components makes use of an action of the centre $Z(G)$ on a~fundamental Weyl alcove~$\Delta$ in
$\mathfrak{t}$, the Lie algebra of a~maximal torus.
The action is described concretely in~\cite{TL} for classical groups and Appendix~\ref{app} records the action for the
two remaining exceptional cases.

Finally, we compute the obstruction to prequantization in Theorem~\ref{thm:obs} in the case $G=\PU(p)$ ($p>2$, prime)
for any number of boundary components~$s$.
We work within the theory of quasi-Hamiltonian group actions with group-valued moment map~\cite{AMM}, where the moduli
space $\mathfrak{M}$ is a~central example.
In quasi-Hamiltonian geometry, quantization is def\/ined as a~certain element of the twisted~$K$-theory
of~$G$~\cite{Khomology}, analogous to $\Spinc$ quantization for Hamiltonian group actions on symplectic manifolds.
In this context, the obstruction to the existence of a~prequantization is a~cohomological obstruction (see
Def\/inition~\ref{preq}).
The obstruction for other cases of non-simply connected structure group does not follow from the approach here (see
Remark~\ref{remark:lament}) and will be considered elsewhere.

\section{Preliminaries}

{\bf Notation.}
Unless otherwise indicated,~$G$ denotes a~compact, simply connected, simple Lie group with Lie algebra $\mathfrak{g}$.
We f\/ix a~maximal torus $T\subset G$ and use the following notation:

\begin{tabular}{l@{\,\,}l} $\mathfrak{t}$ & -- Lie algebra of~$T$;
\\
$\mathfrak{t}^*$ & -- dual of the Lie algebra of~$T$;
\\
$W=N(T)/T$ & -- Weyl group;
\\
$I= \ker \exp_T$ & -- integer lattice;
\\
$P =I^*\subset \mathfrak{t}^*$ & -- (real) weight lattice;
\\
$Q \subset \mathfrak{t}^*$ & -- root lattice;
\\
$Q^\vee\subset \mathfrak{t}$ & -- coroot lattice;
\\
$P^\vee \subset \mathfrak{t}$ & -- coweight lattice.
\end{tabular}

Recall that since~$G$ is simply connected, $I=Q^\vee$.
Moreover, the coroot lattice and weight lattice are dual to each other, as are the root lattice and coweight lattice.
A~choice of simple roots $\alpha_1, \ldots, \alpha_l$ (with $l=\rank(G)$) spanning~$Q$, determines the fundamental
coweights $\lambda_1^\vee, \ldots, \lambda_l^\vee$ spanning $P^\vee$, def\/ined by $\langle\alpha_i,\lambda_j^\vee
\rangle=\delta_{i,j}$.

We let $\langle-,-\rangle$ denote the basic inner product, the invariant inner product on $\mathfrak{g}$ normalized to
make short coroots have length $\sqrt{2}$.
With this inner product, we will often identify $\mathfrak{t} \cong \mathfrak{t}^*$.

Given a~subgroup~$Z$ of the centre $Z(G)$ of~$G$, we shall abuse notation and denote by $q:G\to G/Z$ the resulting
covering(s).

Finally, let $\{e_1, \ldots, e_n\}$ denote the standard basis for ${\mathbb{R}}^n$, equipped with the standard inner
product that will also be denoted with angled brackets $\langle -, - \rangle$.

{\bf Quasi-Hamiltonian group actions.}
We recall some basic def\/initions and facts from~\cite{AMM}.
(For the remainder of this section, we may take~$G$ to be any compact Lie group with invariant inner product
$\langle-,-\rangle$ on $\mathfrak{g}$.) Let $\theta^L$, $\theta^R$ denote the left-invariant, right-invariant
Maurer--Cartan forms on~$G$, and let $\eta= \tfrac{1}{12} \langle \theta^L, [\theta^L,\theta^L] \rangle $ denote the
Cartan 3-form on~$G$.
For a $G$-manifold~$M$, and $\xi \in \mathfrak{g}$, let $\xi^\sharp$ denote the generating vector f\/ield of the action.
The Lie group~$G$ is itself viewed as a $G$-manifold for the conjugation action.

\begin{Definition}[\cite{AMM}]
\label{def:quasi}
A~quasi-Hamiltonian $G$-space is a~triple $(M,\omega,\Phi)$ consisting of a $G$-manifold $M$,
a $G$-invariant 2-form~$\omega$ on~$M$, and an equivariant map $\Phi\colon M \to G$, called the moment map, satisfying:
\begin{enumerate}\itemsep=0pt
\item[i)] $d\omega + \Phi^* \eta = 0 $,
\item[ii)] $\iota_{\xi^\sharp}\omega +\tfrac{1}{2} \Phi^*((\theta^L + \theta^R) \cdot \xi)=0$ for all $\xi \in \mathfrak{g}$,
\item[iii)] at every point $x\in M$, $\ker \omega_x \cap \ker \mathrm{d} \Phi_x = \{0\}$.
\end{enumerate}
\end{Definition}

We will often denote a~quasi-Hamiltonian $G$-space $(M,\omega,\Phi)$ simply by the underlying space~$M$ when~$\omega$
and~$\Phi$ are understood from the context.

The \emph{fusion product} of two quasi-Hamiltonian $G$-spaces $M_j$ with moment maps $\Phi_j:M_j \to G$ ($j=1,2$) is the
product $M_1\times M_2$, with the diagonal $G$-action and moment map $\Phi:M_1 \times M_2 \to G$ given by composing
$\Phi_1 \times \Phi_2$ with multiplication in~$G$.

The \emph{symplectic quotient} of a~quasi-Hamiltonian $G$-space is the symplectic space $M/ /G = \Phi^{-1}(1)/G$, which
is a~symplectic orbifold whenever the group unit $1 \in G$ is a~regular value.
If $1$ is a~singular value, then the symplectic quotient is a~singular symplectic space as def\/ined in~\cite{MS}.

The conjugacy classes $\mathcal{C}\subset G$, with moment map the inclusion into~$G$, are basic examples of
quasi-Hamiltonian $G$-spaces.
Another important example is the \emph{double} $D(G)=G\times G$, equipped with diagonal $G$-action and moment map
$\Phi(g,h)=ghg^{-1}h^{-1}$, the group commutator.
These two families of examples form the building blocks of the moduli space of f\/lat $G$-bundles over a~surface~$\Sigma$
with prescribed boundary holonomies.
(See Section~\ref{sec:comp} for a~sketch of this construction.)

\section{Conjugacy classes invariant under translation \\
by central elements}

This section describes the set of conjugacy classes $\mathcal{D}\subset G$ that are invariant under translation~by
a~subgroup $Z_\mathcal{D}$ of the centre $Z(G)$ of~$G$.
We begin with the following Lemma, which identif\/ies such a~subgroup $Z_\mathcal{D}$ with the fundamental group of
a~conjugacy class in $G/Z$, where $Z\subset Z(G)$.

\begin{Lemma}
\label{lemma:pi1ofC}
Let~$Z$ be a~subgroup of the centre $Z(G)$ of~$G$ and let $\mathcal{C} \subset G/Z$ be a~conjugacy class.
For any conjugacy class $\mathcal{D}\subset G$ covering $\mathcal{C}$, the restriction $q|_{\mathcal{D}}: \mathcal{D}
\to \mathcal{C}$ is the universal covering projection and hence the fundamental group $\pi_1(\mathcal{C}) \cong
Z_\mathcal{D}=\{z \in Z \colon z\mathcal{D} = \mathcal{D} \}$.
\end{Lemma}

\begin{proof}
The inverse image $q^{-1}(\mathcal{C})$ is a~disjoint union of conjugacy classes in~$G$ that cover $\mathcal{C}$.
Since conjugacy classes in a~compact simply connected Lie group are simply connected and $Z_\mathcal{D}$ acts freely on
$\mathcal{D}$, the lemma follows.
\end{proof}

Recall that every element in~$G$ is conjugate to a~unique element $\exp \xi$ in~$T$, where $\xi$ lies in a~f\/ixed (closed)
alcove $\Delta \subset \mathfrak{t}$ of a~Weyl chamber.
Therefore, the set of conjugacy classes in~$G$ is parametrized by~$\Delta$.
Since the $Z(G)$-action commutes with the conjugation action, we obtain an action $Z(G)\times \Delta \to \Delta$.
Next we identify this description of the action of $Z(G)$ on an alcove~$\Delta$ with a~more concrete description of
a~$Z(G)$-action on~$\Delta$ given in~\cite[Section~4.1]{TL}.
(See also~\cite[Section~3.1]{BFM} for a~similar treatment.)

Let $\{\alpha_1, \ldots, \alpha_l \}$ be a~basis of simple roots for $\mathfrak{t}^*$, with highest root
$\tilde{\alpha}=:-\alpha_0$.
Let $\Delta \subset \mathfrak{t}$ be the alcove
\begin{gather*}
\Delta = \{\xi \in \mathfrak{t} \colon \langle \xi,\alpha_j\rangle \geq 0, \, \langle \xi, \tilde{\alpha}\rangle \leq 1
\}.
\end{gather*}

Recall that the exponential map induces an isomorphism $Z(G)\cong P^\vee/Q^\vee$, and that the non-zero elements of the
centre have representatives $\lambda_i^\vee \in P^\vee$ given by minimal dominant coweights.
By~\cite[Lemma~2.3]{TL} the non-zero minimal dominant coweights $\lambda_i^\vee$ are dual to the special roots~$\alpha_i$, which are those roots with coef\/f\/icient~$1$ in the expression $\tilde{\alpha}=\sum m_i\alpha_i$.
In Proposition~4.1.4 of~\cite{TL}, Toledano-Laredo provides a~$Z(G)$-action on~$\Delta$ def\/ined by
\begin{gather*}
z\cdot \xi = w_i \xi + \lambda_i^\vee,
\end{gather*}
where $z=\exp \lambda_i^\vee$, and $w_i \in W$ is a~certain element of the Weyl group.
The element $w_i \in W$ is the unique element that leaves $\Delta \cup \{{\alpha_0} \}$ invariant (i.e.,~induces an
automorphism of the extended Dynkin diagram) and satisf\/ies $w_i({\alpha_0}) = \alpha_i$ (see~\cite[Proposition 4.1.2]{TL}).
The following proposition shows these actions coincide.

\begin{Proposition}
\label{prop:action}
The translation action of $Z(G)$ on~$G$ induces an action $Z(G) \times \Delta \to \Delta$ and is given by the formula
$z\cdot \xi = w_i \xi + \lambda_i^\vee$, where $z=\exp \lambda_i^\vee$ and $w_i$ is the unique element in~$W$ that
leaves $\Delta \cup \{{\alpha_0} \}$ invariant and satisfies $w_i({\alpha_0}) = \alpha_i$.
\end{Proposition}

\begin{proof}
Observe that for any element~$w$ in~$W$, $w \lambda_i^\vee - \lambda_i^\vee \in I= Q^\vee$ since $w\exp\lambda_i^\vee =
\exp \lambda_i^\vee$.
Therefore, $w_i \xi + \lambda_i^\vee = w_i(\xi + \lambda_i^\vee + (w_i^{-1}\lambda_i^\vee - \lambda_i^\vee))$.
In other words, $w_i \xi + \lambda_i^\vee = \hat{w}(\xi+\lambda_i^\vee)$ for some $\hat{w}$ in the af\/f\/ine Weyl group.
Letting $z=\exp\lambda_i^\vee$, this shows that $z \exp \xi = \exp(\xi+\lambda_i^\vee)$ is conjugate to $\exp(w_i \xi
+\lambda^\vee_i)$, which proves the proposition.
\end{proof}

In fact, as the next proposition shows, the automorphism of the Dynkin diagram induced by $w_i$ encodes the resulting
permutation of the vertices of the alcove~$\Delta$.

\begin{Proposition}
\label{prop:dynk}
Let $v_0, \ldots, v_l$ denote the vertices of~$\Delta$ with $v_j$ opposite the facet parallel to $\ker \alpha_j$.
Then $\exp \lambda_i^\vee \cdot v_j = v_k$ whenever $w_i \alpha_j = \alpha_k$, where $w_i$ is as in
Proposition~{\rm \ref{prop:action}}.
\end{Proposition}
\begin{proof}
Let $v_0, \ldots, v_l$ denote the vertices of~$\Delta$, where the vertex $v_j$ is opposite the facet (codimension~$1$
face) parallel to $\ker \alpha_j$.
That is, $v_0=0$ and for $j\neq 0$, $v_j$ satisf\/ies:
\begin{gather*}
\langle \alpha_0,v_j\rangle = -1
\qquad
\text{and}
\qquad
\langle\alpha_r,v_j\rangle = 0
\quad
\text{if and only if}
\quad
0\neq r\neq j.
\end{gather*}
(Hence, for $j\neq 0$ we have $\langle\alpha_j,v_j\rangle = \frac{1}{m_j}$, where $m_j$ is the coef\/f\/icient of $\alpha_j$
in the expression $\tilde\alpha = \sum m_i\alpha_i$.)

Suppose that $w_i\alpha_0 = \alpha_i$ and let $w_i\alpha_j=\alpha_k$ (where~$k$ depends on~$j$).

Consider $\exp \lambda_i^\vee \cdot v_0$.
Since $\langle \alpha_0,w_i v_0+\lambda_i^\vee\rangle = \langle\alpha_0,\lambda_i^\vee\rangle = -1$, and (for $r\neq 0$)
$\langle \alpha_r,w_iv_0 +\lambda_i^\vee\rangle=\langle\alpha_r,\lambda_i^\vee\rangle = \delta_{r,k}$, we have $\exp
\lambda_i^\vee \cdot v_0 = v_i$.

Next, consider $\exp \lambda_i^\vee \cdot v_j=w_i v_j + \lambda_i^\vee$, where $j\neq 0$.
If $k=0$ so that $w_i\alpha_j=\alpha_0$ then $\alpha_j=w_i^{-1} \alpha_0$ is a~special root (i.e.~$m_j=1$) since
$w_i^{-1}=w_j$.
Therefore, $\langle\alpha_0,w_i v_j+\lambda_i^\vee\rangle = \langle w_i^{-1}\alpha_0,v_j\rangle-1= \langle
\alpha_j,v_j\rangle-1=0$.
And if $r\neq 0$,
\begin{gather}
\label{eq:idvert}
\langle\alpha_r,w_i v_j+\lambda_i^\vee\rangle = \langle w_i^{-1} \alpha_r,v_j\rangle + \langle
\alpha_r,\lambda_i^\vee\rangle.
\end{gather}
If $r\neq i$, $w_i^{-1} \alpha_r$ is a~simple root other than $\alpha_j$; therefore, each term above is~$0$.
Moreover, if $r=i$, then the above expression becomes $\langle \alpha_0,v_j \rangle +
\langle\alpha_i,\lambda_i^\vee\rangle = -1+1=0$.
Hence we have $\exp \lambda_i^\vee \cdot v_j = v_0$ whenever $w_i\alpha_j = \alpha_0$.

On the other hand, if $k\neq 0$ so that $w_i \alpha_j = \alpha_k$ is a~simple root, then $\langle \alpha_0,w_i
v_j+\lambda_i^\vee\rangle = \langle w_i^{-1} \alpha_0, v_j \rangle + \langle \alpha_0,\lambda_i^\vee\rangle = 0-1=-1$
since the simple root $w_i^{-1}\alpha_0\neq \alpha_j$.
And if $r\neq 0$, we consider again the expression~\eqref{eq:idvert} and f\/ind (for the same reason as above)
that~\eqref{eq:idvert} is trivial whenever $r\neq k$.
If $r=k$,~\eqref{eq:idvert} becomes $\langle \alpha_k,w_i v_j+ \lambda_i^\vee \rangle = \langle w_i^{-1} \alpha_k,v_j
\rangle + \langle \alpha_k,\lambda_i^\vee \rangle = \langle \alpha_j,v_j \rangle \neq 0$.
Hence we have that $\exp \lambda_i^\vee \cdot v_j = v_k$, as required.
\end{proof}

The $Z(G)$-action on~$\Delta$ is explicitly described in~\cite{TL} for all classical groups.
(In Appendix~\ref{app}, we record the action of the centre on the alcove for the exceptional groups $E_6$ and $E_7$, the
remaining compact simple Lie groups with non-trivial centre.)

{\bf Conjugacy classes in $\boldsymbol{\SU(n)}$.}
We now specialize to the case $G=\SU(n)$ and consider the action of the centre on the alcove.
Identify $\mathfrak{t}\cong \mathfrak{t}^* \subset {\mathbb{R}}^n$ as the subspace $\{x= \sum x_j e_j \colon \sum x_j =0
\}$ and recall that the basic inner product coincides with (the restriction of) the standard inner product on
${\mathbb{R}}^n$.
The roots are the vectors $e_i-e_j$ with $i\neq j$.
Taking the simple roots to be $\alpha_i=e_i-e_{i+1}$ ($i=1,\ldots,n-1$) and the resulting highest root
$\tilde{\alpha}=e_1-e_n$ gives the alcove
\begin{gather*}
\Delta = \{x \in \mathfrak{t} \colon x_1 \geq x_2 \geq \dots \geq x_n, \, x_1-x_n\leq 1 \}.
\end{gather*}
Its vertices are
\begin{gather*}
v_0=0
\qquad
\text{and}
\qquad
v_j=\sum\limits_{i=1}^{j}e_i-\frac{j}{n}\sum\limits_{i=1}^n e_i,
\qquad
j=1, \ldots, n.
\end{gather*}

The centre $Z(\SU(n))\cong {\mathbb{Z}}/n{\mathbb{Z}}$ is generated by ($\exp$ of) the minimal dominant coweight
$\lambda_1^\vee = e_1 -\frac{1}{n}\sum\limits_{i=1}^n e_i$ corresponding to the special root $\alpha_1=e_1-e_2$.
Since the element $w_1$ inducing an automorphism of the extended Dynkin diagram for $\SU(n)$ satisf\/ies
$w_1\alpha_0=\alpha_1$, by Proposition~\ref{prop:dynk} the permutation of the vertices of~$\Delta$ induced by the action
of $\exp \lambda_1^\vee$ is the~$n$-cycle $(v_0   v_1  \cdots   v_{n-1})$ (since $v_j$ is the vertex opposite the
facet parallel to $\ker \alpha_j$).

It follows that the only point in~$\Delta$ f\/ixed by the action of $Z(G)$ is the barycenter
\begin{gather*}
\zeta_* = \frac{1}{n} \sum\limits_{j=0}^{n-1} v_j = \frac{n-1}{2n}e_1 + \frac{n-3}{2n} + \dots + \frac{1-n}{2n}e_n.
\end{gather*}
Hence there is a~unique conjugacy class in $\SU(n)$ that is invariant under translation by the centre~-- namely, matrices
in $\SU(n)$ with eigenvalues $z_1, \ldots, z_n$, the distinct~$n$-th roots of $(-1)^{n+1}$.
As the next proposition shows, however, restricting the action to a~proper subgroup $Z\cong
{\mathbb{Z}}/\nu{\mathbb{Z}}$ ($\nu|n$) of the centre results in larger~$Z$-f\/ixed point sets in~$\Delta$.

\begin{Proposition}
Let $n=\nu m$ and consider the subgroup ${\mathbb{Z}}/\nu{\mathbb{Z}} \subset {\mathbb{Z}}/n{\mathbb{Z}} \cong
Z(\SU(n))$.
The ${\mathbb{Z}}/\nu{\mathbb{Z}}$-fixed points in the alcove~$\Delta$ for $\SU(n)$ consist of the convex hull of the
barycenters of the faces spanned by the orbits of the vertices $v_0,\ldots,v_{m-1}$ of~$\Delta$.
\end{Proposition}

\begin{proof}
Write $x=\sum t_i v_i$ in~$\Delta$ in barycentric coordinates (with $t_i\geq 0$ and $\sum t_i=1$).
Then a~generator of ${\mathbb{Z}}/\nu{\mathbb{Z}}$ sends~$x$ to $\sum t_i' v_i$, with $t_i'=t_{i-m \mod n}$.
Therefore~$x$ is f\/ixed if and only if $t_i=t_{i-m \mod n}$, and in this case we may write,
\begin{gather*}
x=t_0 \sum\limits_{j=0}^{\nu-1} v_{jm} + t_1 \sum\limits_{j=0}^{\nu-1} v_{1+jm} + \dots + t_{m-1}
\sum\limits_{j=0}^{\nu-1} v_{m-1+jm}
\\
\phantom{x}{}
 =\nu t_0 \frac{1}{\nu}\sum\limits_{j=0}^{\nu-1} v_{jm} + \nu t_1 \frac{1}{\nu} \sum\limits_{j=0}^{\nu-1} v_{1+jm} +
\dots + \nu t_{m-1} \frac{1}{\nu}\sum\limits_{j=0}^{\nu-1} v_{m-1+jm},
\end{gather*}
which exhibits a~f\/ixed point in the desired form.
\end{proof}

\begin{figure}[h] \centering
\includegraphics[scale=0.95]{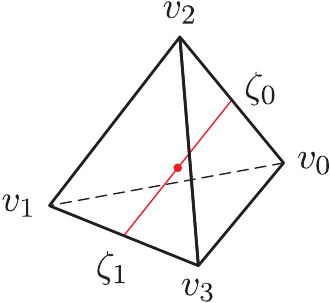}

\caption{Alcove for $\SU(4)$.
The indicated line segment through the barycenter parametrizes the set of conjugacy classes invariant under translation
by ${\mathbb{Z}}/2{\mathbb{Z}} \subset Z(\SU(4))$.}
\label{fig:su4}
\end{figure}

To illustrate, consider the subgroup $Z \cong {\mathbb{Z}}/2{\mathbb{Z}}$ of the centre $Z(\SU(4)) \cong
{\mathbb{Z}}/4{\mathbb{Z}}$, which acts by transposing the vertices $v_0 \leftrightarrow v_2$ and $v_1 \leftrightarrow
v_3$.
The barycenters $\zeta_0$, $\zeta_1$ of the edges~$\overline{v_0v_2}$ and~$\overline{v_1v_3}$, respectively, are f\/ixed
and thus the~$Z$-f\/ixed points are those on the line segment joining~$\zeta_0$ and~$\zeta_1$
(see Fig.~\ref{fig:su4}).

\section{Components of the moduli space with markings}\label{sec:comp}

In this section we recall the quasi-Hamiltonian description of the moduli space of f\/lat bundles over a~compact
orientable surface with prescribed boundary holonomies.
We refer to the original article~\cite{AMM} for the   details regarding the construction sketched below.

Let~$\Sigma$ be a~compact, oriented surface of genus~$h$ with~$s$ boundary components.
For conjugacy classes $\mathcal{C}_1, \ldots, \mathcal{C}_s$ in $G/Z$, let $M_{G/Z}(\Sigma; \mathcal{C}_1, \ldots,
\mathcal{C}_s)$ be the moduli space of f\/lat $G/Z$-bundles over~$\Sigma$ with prescribed boundary holonomies lying in the
conjugacy classes $\mathcal{C}_j$ ($j=1, \ldots, s$).
Points in $M_{G/Z}(\Sigma; \mathcal{C}_1, \ldots, \mathcal{C}_s)$ are (gauge equivalence classes of) principal
$G/Z$-bundles over~$\Sigma$ equipped with a~f\/lat connection whose holonomy around the~$j$-th boundary component lies in
the conjugacy class $\mathcal{C}_j$.
This moduli space is an important example in the theory of quasi-Hamiltonian group actions, where it is cast
a~symplectic quotient of a~fusion product,
\begin{gather}
\label{eq:moduli}
M_{G/Z}(\Sigma; \mathcal{C}_1, \ldots, \mathcal{C}_s) = \big(D(G/Z)^h\times \mathcal{C}_1 \times \dots \times \mathcal{C}_s\big)
/ /(G/Z),
\end{gather}
which may have several connected components if~$Z$ is non-trivial.
Extending the discussion in~\cite[Section 2.3]{KM}, we describe the connected components of~\eqref{eq:moduli} as
symplectic quotients of an auxiliary quasi-Hamiltonian $G$-space.

As in~\cite[Section 2.2]{KM}, given a~quasi-Hamiltonian $G/Z$-space~$N$ with group-valued moment map $\Phi:N\to G/Z$,
let $\check{N}$ be the f\/ibre product def\/ined by the Cartesian square,
\begin{gather}
\label{diag:Ncheck}
\begin{split}
& \xymatrix{
\check{N} \ar[d] \ar[r]^{\check{\Phi}} & G \ar[d]^q \\
N \ar[r]^{\Phi} & G/Z
}
\end{split}
\end{gather}
Then $\check{N}$ is naturally a~quasi-Hamiltonian $G$-space with moment map $\check{\Phi}$.
The following proposition from~\cite{KM} and its Corollary summarize some properties of this construction.

\begin{Proposition}[\protect{\cite[Proposition~2.2]{KM}}]
\label{prop:pb}
Let $\check{N}$ be the fibre product defined by~\eqref{diag:Ncheck}, where $\Phi:N\to G/Z$ is a~group-valued moment map.

\begin{enumerate}\itemsep=0pt
\item[$i)$]
We have a~canonical identification of symplectic quotients $\check{N}/ / G \cong N/ / (G/Z)$.
\item[$ii)$]
For a~fusion product $N={N}_1\times \dots \times {N}_r$ of quasi-Hamiltonian ${G/Z}$-spaces, the space $\check{N}$ is
a~quotient of $\check{N}_1\times\dots \times \check{N}_r$ by the group $\big\{(c_1,\ldots,c_r)\in Z^r \mid \prod\limits_{j=1}^r c_j=e\big\}$.
\item[$iii)$]
If $\Phi\colon N\to G/Z$ lifts to a~moment map $\Phi'\colon N\to {G}$, thus turning~$N$ into a~quasi-Hamiltonian
${G}$-space   then $ \check{N}=N\times Z$.
\end{enumerate}
\end{Proposition}

\begin{Corollary}
\label{cor:comps}
Let $\check{N}$ be the fibre product defined by~\eqref{diag:Ncheck}, where $\Phi:N\to G/Z$ is a~group-valued moment
map, and write $\check{N}=\bigsqcup X_j$ as a~union of its connected components.
Then the components of $N/ /(G/Z)$ can be identified with the symplectic quotients $X_j/ /G$.
\end{Corollary}

\begin{proof}
The restrictions $\check{\Phi}_j=\check{\Phi}|_{X_j}$ are $G$-valued moment maps whose f\/ibres are connected by \cite[Theorem~7.2]{AMM}.
Since $\check{\Phi}^{-1}(e) = \bigsqcup \check{\Phi}_j^{-1}(e)$, it follows that $\check{N}/ /G = \check{\Phi}^{-1}(e)/G
= \bigsqcup \check{\Phi}_j^{-1}(e)/G=\bigsqcup X_j/ /G$.
The result follows from Proposition~\ref{prop:pb}(i). 
\end{proof}

Hence to identify the components of~\eqref{eq:moduli}, it suf\/f\/ices to identify the components of $\check{N} / /G$, where
$N= D(G/Z)^h\times\mathcal{C}_1 \times \dots \times \mathcal{C}_s $~-- namely, $X_j/ /G$, where $X_j$ ranges over the
components of~$\check{N}$.
In particular, we may view the moduli space~\eqref{eq:moduli} as a~union of symplectic quotients of
quasi-Hamiltonian $G$-spaces (as opposed to $G/Z$-spaces), which will be very important for the approach taken in
Section~\ref{sec:main}.

With this in mind, choose conjugacy classes $\mathcal{D}_j \subset G$ covering $\mathcal{C}_j$ ($j=1,\ldots, s$) and let
\begin{gather*}
\widetilde{N}= D(G)^h\times\mathcal{D}_1 \times \cdots \times \mathcal{D}_s.
\end{gather*}
Let
\begin{gather}
\label{eq:gamma}
\Gamma = \Big\{(\gamma_1, \ldots, \gamma_{s}) \in Z_{\mathcal{D}_1} \times \dots \times Z_{\mathcal{D}_s} \colon \prod
\gamma_j =1 \Big\} \subset Z^s
\end{gather}
(cf.~Lemma~\ref{lemma:pi1ofC}).
We show next that the components of $\check{N}$ are all homeomorphic to $\widetilde{N}/(Z^{2h}\times \Gamma)$
(generalizing the decomposition appearing in~\cite[Lemma~2.3]{KM} for $G=\SO(3)$).

\begin{Proposition}
\label{prop:comps}
Let $N=D(G/Z)^h\times \mathcal{C}_1 \times \dots \times \mathcal{C}_s $ for conjugacy classes $\mathcal{C}_j \subset
G/Z$ $(j=1,\ldots,s)$ and let $\check{N}$ be the fibre product defined by~\eqref{diag:Ncheck}.
Then $\check{N}$ may be written as a~union of its connected components,
\begin{gather*}
\check{N} \cong \bigsqcup_{Z/(Z_{\mathcal{D}_1} \dots Z_{\mathcal{D}_s})}D(G/Z)^h\times (\mathcal{D}_1 \times \cdots
\times \mathcal{D}_s)/{\Gamma},
\end{gather*}
where $\mathcal{D}_j \subset G$ are conjugacy classes covering $\mathcal{C}_j$ $(j=1,\ldots,s)$ and~$\Gamma$ is as
in~\eqref{eq:gamma}.
\end{Proposition}
\begin{proof}
This is a~straightforward application of the properties~(ii) 
and~(iii) 
listed in Proposition~\ref{prop:pb}.
By property~(iii), $\check{D}(G/Z)^h = D(G/Z)^h \times Z$, and by Lemma~\ref{lemma:pi1ofC},
$\check{\mathcal{C}}_j = \mathcal{D}_j \times Z/Z_{\mathcal{D}_j}$.
Therefore, by property~(ii),
\begin{gather*}
\check{N} \cong D(G/Z)^h \times (Z\times \mathcal{D}_1 \times Z/Z_{\mathcal{D}_1} \times \cdots \times \mathcal{D}_s
\times Z/Z_{\mathcal{D}_s})/\Lambda,
\end{gather*}
where $\Lambda={\{(c_0, \ldots, c_{s}) \in Z^{s+1} \colon c_0 \cdots c_{s}=1\}}$.
Since
\begin{gather*}
(Z\times \mathcal{D}_1 \times Z/Z_{\mathcal{D}_1} \times \dots \times \mathcal{D}_s \times Z/Z_{\mathcal{D}_s})/\Lambda
\cong (Z\times \mathcal{D}_1 \times \dots \times \mathcal{D}_s)/\Gamma',
\end{gather*}
where $\Gamma'=\{(\gamma_0, \ldots, \gamma_{s}) \in Z\times Z_{\mathcal{D}_1} \times \dots \times Z_{\mathcal{D}_s}
\colon \prod \gamma_j =1\}$, we see that the components of $\check{N}$ are in bijection with $Z/(Z_{\mathcal{D}_1} \cdots
Z_{\mathcal{D}_s})$.

Consider the component corresponding to $\bar{z} \in Z/(Z_{\mathcal{D}_1}\cdots Z_{\mathcal{D}_s})$ in which each point
is of the form $(\vec{g},[(z,x_1, \ldots, x_s)]_{\Gamma'})$, where $[\;]_{\Gamma'}$
denotes a~$\Gamma'$-orbit.
(Note that there is always a~representative of this form with~$z$ in the f\/irst coordinate.) This component is
homeomorphic to $D(G/Z)^h \times (\mathcal{D}_1 \times \cdots \times \mathcal{D}_s)/\Gamma$ by the map $(\vec{g},[(x_1, \ldots,
x_s)]_{\Gamma}) \mapsto (\vec{g},[(z,x_1, \ldots, x_s)]_{\Gamma'})$.
\end{proof}

\begin{Remark}
The decomposition in Proposition~\ref{prop:comps} is consistent with Theorem~14 in Ho--Liu's work~\cite{HL} on the
connected components of the moduli space for any compact connected Lie group~$G$.
\end{Remark}

For the case $G/Z=\SU(p)/({\mathbb{Z}}/p{\mathbb{Z}})=\PU(p)$, where~$p$ is prime, the decomposition above simplif\/ies.
In particular, there is only one conjugacy class $\mathcal{D}_*=\SU(p) \cdot \exp \zeta_*$, corresponding to the
barycenter $\zeta_* \in \Delta$, invariant under the action of the centre.
Let $\mathcal{C}_*=q(\mathcal{D}_*)$ be the corresponding conjugacy class in $\PU(p)$.
Therefore, we obtain the following Corollary (cf.~\cite[Lemma 2.3]{KM}).

\begin{Corollary}
\label{cor:PUpdecomp}
Let~$p$ be prime and let $N= D(\PU(p))^h \times \mathcal{C}_1 \times \dots \times \mathcal{C}_s $ for conjugacy classes
$\mathcal{C}_j \subset \PU(p)$ $(j=1,\ldots,s)$ and let $\check{N}$ be the fibre product defined by~\eqref{diag:Ncheck}.
Then,
\begin{gather*}
\check{N} \cong
\begin{cases}
D(\PU(p))^h \times (\mathcal{D}_1 \times \cdots \times \mathcal{D}_s)/{\Gamma} & \text{if}\quad\exists\, j \colon \mathcal{C}_j=\mathcal{C}_*,
\\
D(\PU(p))^h \times {\mathcal{D}_1 \times \cdots \times \mathcal{D}_s} \times Z & \text{otherwise},
\end{cases}
\end{gather*}
where $\mathcal{D}_j \subset \SU(p)$ are conjugacy classes covering $\mathcal{C}_j$ $(j=1,\ldots,s)$ and~$\Gamma$ is as
in~\eqref{eq:gamma}.
\end{Corollary}

In particular, if (after re-labelling) $\mathcal{C}_j=\mathcal{C}_*$ for all $j\leq r$ ($r>0$), then we obtain
\begin{gather}
\label{eq:splits4PUp}
\check{N} \cong D(\PU(p))^h \times (\mathcal{D}_*)^r/{\Gamma} \times \mathcal{D}_{r+1} \times \cdots \times
\mathcal{D}_s,
\end{gather}
where, in this case, $\Gamma=\{(\gamma_1, \ldots, \gamma_r) \in Z^r \colon \prod \gamma_j=1\}$.

\section{Obstruction to prequantization}\label{sec:main}

\subsection{Prequantization for quasi-Hamiltonian group actions}\label{sec:preq}

We recall some def\/initions and properties regarding prequantization of quasi-Hamiltonian group actions.
Recall that the Cartan $3$-form $\eta \in \Omega^3(G)$ is integral~-- in fact, $[\eta] \in H^3(G;{\mathbb{R}})$ is the
image of a~generator $x\in H^3(G;{\mathbb{Z}})\cong {\mathbb{Z}}$ under the coef\/f\/icient homomorphism induced~by
${\mathbb{Z}} \to {\mathbb{R}}$.
Condition~(i) 
in Def\/inition~\ref{def:quasi} says that the pair $(\omega,\eta)$ def\/ines a~relative
cocycle in $\Omega^3(\Phi)$, the algebraic mapping cone of the pull-back map $\Phi^*\colon \Omega^*(G)\to\Omega^*(M)$,
and hence a~cohomology class $[(\omega,\eta)] \in H^3(\Phi;{\mathbb{R}})$.
(See~\cite[Chapter~I, Section~6]{BottTu} for the def\/inition of relative cohomology.)
\begin{Definition}[\cite{JGP, Khomology}]\label{preq}
Let $k\in {\mathbb{N}}$.
A~\emph{level~$k$ prequantization} of a~quasi-Hamiltonian $G$-space $(M,\omega,\Phi)$ is an integral lift $\alpha\in
H^3(\Phi;{\mathbb{Z}})$ of the class $k[(\omega,\eta)]\in H^3(\Phi;{\mathbb{R}})$.
\end{Definition}

\begin{Remark}
The def\/inition of prequantization in Def\/inition~\ref{preq} uses the assumption in this paper that~$G$ is simply
connected.
The general def\/inition of prequantization~\cite[Def\/inition~3.2]{Khomology} (with~$G$ semi-simple and compact) requires
an integral lift in $H^3_G(\Phi;{\mathbb{Z}})$ of an equivariant extension of the class $k[(\omega,\eta)]$.
When~$G$ is simply connected,~\cite[Proposition~3.5]{JGP} shows that the def\/inition above is equivalent.
Our main goal is to apply the quasi-Hamiltonian viewpoint on prequantization to the moduli space of f\/lat bundles with
prescribed holonomies; therefore, by Corollary~\ref{cor:comps} it suf\/f\/ices to work on each component of the
quasi-Hamiltonian $G$-space $\check{N}$ in Proposition~\ref{prop:comps} using Def\/inition~\ref{preq}.
\end{Remark}

We list some basic properties of level~$k$ prequantizations that we shall encounter.
\begin{enumerate}\itemsep=0pt
\item[(a)] If $M_1$ and $M_2$ are pre-quantized quasi-Hamiltonian $G$-spaces at level~$k$, then their fusion product
$M_1 \times M_2$ inherits a~prequantization at level~$k$.
Conversely, a~prequantization of the product induces prequantizations of the factors.
See~\cite[Proposition~3.8]{JGP}.
\item[(b)] A~level~$k$ prequantization of~$M$ induces a~prequantization of the symplectic quotient $M/  / G$, equipped
with the~$k$-th multiple of the symplectic form.
\item[(c)] The long exact sequence in relative cohomology gives a~necessary condition $k\Phi^*(x)=0$ for the existence
of a~level~$k$-prequantization.
If $H^2(M;{\mathbb{R}})=0$, $k\Phi^*(x)=0$ is also suf\/f\/icient~\cite[Proposition~4.2]{JGP} to conclude
a~level~$k$-prequantization exists.
\end{enumerate}

The following examples relate to the moduli space of f\/lat bundles with prescribed boundary holonomies.

\begin{Example}
\label{eg:double}
The double $D(G)=G\times G$ with moment map $\Phi: D(G) \to G$ equal to the group commutator admits a~prequantization at
all levels $k\in {\mathbb{N}}$.
For non-simply connected groups, the double $D(G/Z)$ with moment map $\Phi: D(G/Z) \to G$ the canonical lift of the
group commutator admits a~level~$k$-prequantization if and only if~$k$ is a~multiple of $l_0 \in {\mathbb{N}}$, where
$l_0$ is a~positive integer depending on $G/Z$ computed for all compact simple Lie groups in~\cite{JGP}.
For $G/Z = \PU(n)$, $l_0=n$.
\end{Example}

\begin{Example}
Conjugacy classes $\mathcal{D} \subset G$ admitting a~level~$k$-prequantization are those $\mathcal{D}=G\cdot \exp \xi$
($\xi \in \Delta$) with $(k\xi)^\flat \in P$~\cite{BasicGerbe}, where $(k\xi)^\flat = \langle k\xi,- \rangle$ (i.e., a~\emph{level~$k$ weight}).
For simply laced groups (such as $G=\SU(n)$), under the identif\/ication $\mathfrak{t} \cong \mathfrak{t}^*$, $P^\vee
\cong P$.
Therefore, in this case, $\mathcal{D}$ admits a~level~$k$-prequantization if and only if $k\xi \in P^\vee$.
Since $\exp^{-1} Z(G) = P^\vee$, we see that $\mathcal{D}$ admits a~level~$k$-prequantization if and only if $g^k \in
Z(G)$ for all $g \in \mathcal{D}$.
(So in particular if~$k$ is a~multiple of the \emph{order of $\mathcal{D}$}~\cite[Def\/inition 5.76]{BL}, then
$\mathcal{D}$ admits a~level~$k$ prequantization.)
\end{Example}

\subsection{Quasi-Hamiltonian prequantization and symplectic quotients}\label{pqbundle}

To provide some context, we further elaborate on property (b) following Def\/inition~\ref{preq} since we view the moduli
space of f\/lat bundles as a~symplectic quotient~\eqref{eq:moduli} in quasi-Hamiltonian geometry.
By~\cite[Proposition~3.6]{JGP}, a~level~$k$ prequantization of a~quasi-Hamiltonian $G$-space $(M,\omega,\Phi)$ gives
an integral lift of the equivariant cohomology class $[j^*\omega] \in H^2_G(\Phi^{-1}(1)\,;{\mathbb{R}})$, where
$j:\Phi^{-1}(1) \to M$ denotes inclusion.
Hence, there is a $G$-equivariant line bundle $L\to \Phi^{-1}(1)$ with connection of curvature $j^*\omega$.
If~$1$ is a~regular value, the symplectic quotient $M/ /G= \Phi^{-1}(1)/G$ is a~symplectic orbifold~\cite{AMM} and
the $G$-equivariant line bundle over the level set descends to a~prequantum orbifold line bundle $L/G \to M/ /G$.
(See~\cite[Example 2.29]{Adem} for a~discussion of orbifold vector bundles in this context.)

\begin{Remark}
The orbifold line bundle $L/G\to M/ /G$ need not be an ordinary line bundle over the underlying topological space $M/
/G$ (i.e., the coarse moduli space of the orbifold).
(For orbifolds that arise as quotients $X/G$ of a~smooth, proper, locally free action of a~Lie group~$G$ on a~smooth
manifold~$X$, this distinction is apparent from the observation that $H_G^2(X\,;{\mathbb{Z}})$ is not necessarily
isomorphic to $H^2(X/G\,;{\mathbb{Z}})$.) Some works in the literature require a~prequantization to be an ordinary line
bundle, and hence obtain a~further obstruction to the existence of a~prequantization (e.g.~\cite[Theorem 4.2]{Ko},
\cite[Theorems~4.12 and~6.1]{DW}, \cite[Lemme~3.2]{Pa}).
We gratefully acknowledge the referee's comments that led to this important clarif\/ication.
\end{Remark}

\subsection[The obstruction to prequantization for the moduli space of $\PU(p)$ bundles, $p$ prime]{The obstruction
to prequantization for the moduli space\\ of $\boldsymbol{\PU(p)}$ bundles, $\boldsymbol{p}$ prime}

Let~$p$ be an odd prime.
In this section we obtain the obstruction to prequantization for the quasi-Hamiltonian $\SU(p)$-space $\check{N}$, where
$N= D(\PU(p))^h \times \mathcal{C}_1 \times \dots \times \mathcal{C}_s $ for conjugacy classes $\mathcal{C}_j \subset
\PU(p)$ ($j=1,\ldots,s$).
Let $M\subseteq \check{N}$ be a~connected component (by Corollary~\ref{cor:PUpdecomp}),
\begin{gather*}
M=D(\PU(p))^h \times \left(\mathcal{D}_1 \times \dots \times \mathcal{D}_s\right)/{\Gamma},
\end{gather*}
where~$\Gamma$ is as in~\eqref{eq:gamma}.
As we shall see in the proof of Theorem~\ref{thm:obs}, we will f\/ind property (a) in Section~\ref{sec:preq} very useful
in order to proceed `factor by factor',
using the decomposition~\eqref{eq:splits4PUp}.

To begin, we establish the following proposition which allows us to use property~(c) in Section~\ref{sec:preq} to
compute the obstruction to prequantization for the factor $(\mathcal{D}_*)^r/\Gamma$ in~\eqref{eq:splits4PUp}.

\begin{Proposition}
\label{prop:nofreeH2}
Let $\mathcal{D}_* \subset \SU(p)$ denote the conjugacy class of the barycenter $\zeta_* $ of the alcove~$\Delta$ and
let $\Gamma = \{(\gamma_1, \ldots, \gamma_r) \in Z^r \colon \prod \gamma_j =1 \}$ with $r>1$.
Then $H^2((\mathcal{D}_*)^r/\Gamma;{\mathbb{R}})=0$.
\end{Proposition}
\begin{proof}
Since $(\mathcal{D}_*)^r \to (\mathcal{D}_*)/\Gamma$ is a~covering projection,
$H^2((\mathcal{D}_*)^r/\Gamma;{\mathbb{R}}) \cong H^2 ((\mathcal{D}_*)^r; {\mathbb{R}})^\Gamma$.
By the K\"unneth Theorem, $H^2((\mathcal{D}_*)^r;{\mathbb{R}}) \cong \bigoplus H^2(\mathcal{D}_*,{\mathbb{R}})$.
Since the~$\Gamma$-action factors through $Z^m$, $H^2 ((\mathcal{D}_*)^r; {\mathbb{R}})^\Gamma = \bigoplus
H^2(\mathcal{D}_*;{\mathbb{R}})^Z$.

Recall that since $\zeta_*$ lies in the interior of the alcove, the centralizer $\SU(p)_{\exp \zeta_*} = T$ and hence
$\mathcal{D}_* \cong \SU(p)/T$.
Moreover, we have $H^*(\mathcal{D}_*;{\mathbb{R}})\cong {\mathbb{R}}[t_1, \ldots, t_p]/(\sigma_1, \ldots, \sigma_p)$,
where $\sigma_i$'s are the elementary symmetric polynomials.
In particular, we may write
\begin{gather*}
H^2(\mathcal{D}_*;{\mathbb{R}}) \cong ({\mathbb{R}} t_1 \oplus \dots \oplus {\mathbb{R}} t_p)/(t_1+ \dots + t_p=0).
\end{gather*}
The~$Z$-action on $\mathcal{D}_*$ corresponds to an action on $\SU(p)/T$ by a~cyclic subgroup of the Weyl group (e.g.,
see the proof of Proposition~\ref{prop:action}).
Since the Weyl group (i.e., symmetric group $\Sigma_p$) acts by permuting the $t_i$,~$Z$ acts by a~$p$-cycle on the $t_i$.
Therefore, $H^2(\mathcal{D}_*;{\mathbb{R}})^Z=0$, which establishes the result.
\end{proof}

\begin{Remark}
\label{remark:lament}
The analogue of Proposition~\ref{prop:nofreeH2} for the factors $(\mathcal{D}_1 \times \dots \times
\mathcal{D}_s)/\Gamma$ that appear in the decomposition in Proposition~\ref{prop:comps} need not hold when considering
other non-simply connected structure groups $G/Z$.
\end{Remark}

\begin{Theorem}
\label{thm:obs}
The quasi-Hamiltonian $\SU(p)$-space $M=D(\PU(p))^h \times (\mathcal{D}_1 \times \dots \times
\mathcal{D}_s)/{\Gamma} $ admits a~level~$k$-prequantization if and only if the following conditions are
satisfied:
\begin{enumerate}\itemsep=0pt
\item[$i)$]
\label{thm:obs-double}
if $h\geq 1$, then $k \in p{\mathbb{N}}$;
\item[$ii)$]
\label{thm:obs-conj}
$g^k \in Z(\SU(p))$ for every $\displaystyle g \in \mathcal{D}_1 \cup \cdots \cup \mathcal{D}_s$.
\end{enumerate}
Moreover, if in addition the identity matrix $1 \in \SU(p)$ is a~regular value of the restriction of the group-valued
moment map $\check{\Phi}:M\to \SU(p)$, the prequantization descends to a~prequantization of the corresponding component
of the moduli space $\mathfrak{M}=M_{\PU(p)}(\Sigma; \mathcal{C}_1, \ldots, \mathcal{C}_s)$, where
$\mathcal{C}_j=q(\mathcal{D}_j) \subset \PU(p)$.
\end{Theorem}
\begin{proof}
By property (a) in Section~\ref{sec:preq},~$M$ admits a~level~$k$-prequantization if and only if each factor does.
Since $D(\PU(p))$ admits a~level~$k$-prequantization if and only if condition~(i) 
is satisf\/ied (see Example~\ref{eg:double}), we may assume from now on $h=0$.

We f\/irst verify the necessity of condition~(ii). 
A~prequantization of $M=(\mathcal{D}_1 \times \dots \times \mathcal{D}_s) /\Gamma$ induces a~prequantization of its
universal cover $\tilde{M}=\mathcal{D}_1 \times \dots \times \mathcal{D}_s$, and hence each $\mathcal{D}_j$ must admit
a~prequantization, which is equivalent to condition~(ii). 

Next we verify that condition~(ii) 
is suf\/f\/icient for a~level~$k$-prequantization of~$M$ (with $h=0$).
As in the decomposition~\eqref{eq:splits4PUp}, write (possibly after re-labelling)
\begin{gather*}
M=\underbrace{(\mathcal{D}_* \times \cdots \times \mathcal{D}_*}_{r~\text{factors}})/\Gamma \times \mathcal{D}_{r+1}
\times \cdots \times \mathcal{D}_s
\end{gather*}
Using property (a) in Section~\ref{sec:preq} again, it suf\/f\/ices to consider the case $1<r=s$.
(Note that if $s=r=1$,~$\Gamma$ is trivial.) In this case, condition~(ii) 
is simply that $\mathcal{D}_*$ admit a~level~$k$ prequantization.
Since $\mathcal{D}_*$ consists of matrices in $\SU(p)$ conjugate to
\begin{gather*}
\exp \zeta_* = \mathrm{diag}\big(\exp\big(\tfrac{p-1}{p}\pi\sqrt{-1}\big), \exp\big(\tfrac{p-3}{p}\pi\sqrt{-1}\big), \dots,
\exp\big(\tfrac{1-p}{p}\pi\sqrt{-1}\big)\big),
\end{gather*}
$\mathcal{D}_*$ admits a~level~$k$-prequantization if and only if $(\exp \zeta_*)^k$ is a~scalar matrix; if and only
if~$k$ is a~multiple of~$p$.
By property~(c) in Section~\ref{sec:preq} and Proposition~\ref{prop:nofreeH2}, it suf\/f\/ices to show that
$p\cdot\check{\Phi}^*x =0$, where $\check{\Phi}:M\to \SU(p)$ is the group-valued moment map.

By Corollary~7.6 in~\cite{AM}, $\mathsf{h}^\vee \check{\Phi}^*x=W_3(M)$, the third integral Stiefel--Whitney class, where
$\mathsf{h}^\vee$ denotes the dual Coxeter number.
Recall that $W_3(M)=\beta w_2(M)$, where $\beta: H^2(M;{\mathbb{Z}}/2{\mathbb{Z}}) \to H^3(M;{\mathbb{Z}})$ is the
(integral) Bockstein homomorphism and $w_2(M)$ is the second Stiefel--Whitney class.
Since~$\Gamma$ has odd order, $H^2(M;{\mathbb{Z}}/2{\mathbb{Z}}) \cong
H^2((\mathcal{D}_*)^r;{\mathbb{Z}}/2{\mathbb{Z}})^\Gamma$, which is trivial (by an argument similar to the proof of
Proposition~\ref{prop:nofreeH2}).
Since $\mathsf{h}^\vee=p$, we are done.

As discussed in Section~\ref{pqbundle}, a~quasi-Hamiltonian prequantization descends to a~prequantization of the
symplectic quotient $\mathfrak{M}$.
\end{proof}

\appendix

\section{The action of the centre on the alcove\\
of exceptional Lie groups}\label{app}

Below we record the action of the centre $Z(G)$ on an alcove for the exceptional Lie groups $G=E_6$ and $G=E_7$.
(The action for classical groups appears in~\cite{TL}.)

The vertices of the alcove were obtained using \texttt{polymake}~\cite{polymake}, which outputs the vertices of
a~polytope presented as an intersection of half-spaces.
The relevant Weyl group element from Proposition~\ref{prop:action}~-- one which gives an automorphism of the extended
Dynkin diagram~-- was found with the help of John Stembridge's \texttt{coxeter-weyl} package for Maple~\cite{St};
a~direct calculation then shows that this element has the desired properties in Proposition~\ref{prop:action}.

Let $\{e_1, \ldots, e_8 \}$ denote the standard basis in ${\mathbb{R}}^8$, equipped with the usual inner product.
Given a~vector~$\alpha$ in ${\mathbb{R}}^8$, $\refl{\alpha}:{\mathbb{R}}^8\to {\mathbb{R}}^8$ denotes ref\/lection in the
subspace orthogonal to~$\alpha$,
\begin{gather*}
\refl{\alpha}(v)=v-\frac{2\langle \alpha,v \rangle}{\langle \alpha,\alpha \rangle} \alpha.
\end{gather*}
The notation used below is consistent with that found in~\cite[Planches V-VI]{Bourbaki}.

{\bf $\boldsymbol{G=E_6}$.}
Let $\mathfrak{t} \cong \mathfrak{t}^*\cong \{(x_1, \ldots, x_8) \in {\mathbb{R}}^8 \colon x_6=x_7=-x_8 \}$.
The simple roots $\alpha_1, \ldots, \alpha_6$ and highest root $\tilde\alpha$ determine the half-spaces whose
intersection is the alcove $\Delta \subset \mathfrak{t}$.
The vertices of~$\Delta$ (opposite the facets parallel to the corresponding root hyperplanes) are given in
Table~\ref{table:alcE6}.

\begin{table}[h] \centering\renewcommand{\arraystretch}{1.3}
\caption{Alcove data for $E_6$.}\label{table:alcE6}
\vspace{1mm}

\begin{tabular}{ll}
\hline
Simple or dominant root & Opposite vertex
\\
\hline
$\alpha_1=\big(\frac{1}{2},-\frac{1}{2},-\frac{1}{2},-\frac{1}{2},-\frac{1}{2},-\frac{1}{2},-\frac{1}{2},\frac{1}{2}\big)$ &
$v_1=\big(0,0,0,0,0,-\frac{2}{3},-\frac{2}{3},\frac{2}{3}\big)$
\\
$\alpha_2=(1,1,0,0,0,0,0,0)$ &
$v_2=\big(\frac{1}{4},\frac{1}{4},\frac{1}{4},\frac{1}{4},\frac{1}{4},-\frac{1}{4},-\frac{1}{4},\frac{1}{4}\big)$
\\
$\alpha_3=(-1,1,0,0,0,0,0,0)$ &
$v_3=\big({-}\frac{1}{4},\frac{1}{4},\frac{1}{4},\frac{1}{4},\frac{1}{4},-\frac{5}{12},-\frac{5}{12},\frac{5}{12}\big)$
\\
$\alpha_4=(0,-1,1,0,0,0,0,0)$ & $v_4=\big(0,0,\frac{1}{3},\frac{1}{3},\frac{1}{3},-\frac{1}{3},-\frac{1}{3},\frac{1}{3}\big)$
\\
$\alpha_5=(0,0,-1,1,0,0,0,0)$ & $v_5=\big(0,0,0,\frac{1}{2},\frac{1}{2},-\frac{1}{3},-\frac{1}{3},\frac{1}{3}\big)$
\\
$\alpha_6=(0,0,0,-1,1,0,0,0)$ & $v_6=\big(0,0,0,0,1,-\frac{1}{3},-\frac{1}{3},\frac{1}{3}\big)$
\\
$\tilde{\alpha}=\big(\frac{1}{2},\frac{1}{2},\frac{1}{2},\frac{1}{2},\frac{1}{2},-\frac{1}{2},-\frac{1}{2},\frac{1}{2}\big)$ &
$v_0=0$
\\
\hline
\end{tabular}
\end{table}

The non-zero elements of the centre $Z(E_6) \cong {\mathbb{Z}}/3{\mathbb{Z}}$ are given by ($\exp$ of) the minimal
dominant coweights $\lambda_1^\vee = \frac{2}{3}(e_8-e_7-e_6)$ and $\lambda_6^\vee=e_5+ \frac{1}{3}(e_8-e_7-e_6)$.
The corresponding elements $w_1$ and $w_6$ of the Weyl group (as in Proposition~\ref{prop:action}), inducing
automorphisms of the extended Dynkin diagram are:
\begin{gather*}
w_1=\refl{\alpha_1}\refl{\alpha_3}\refl{\alpha_4}\refl{\alpha_2}\refl{\alpha_5}\refl{\alpha_4}\refl{\alpha_3}\refl{\alpha_1}
    \refl{\alpha_6}\refl{\alpha_5}\refl{\alpha_4}\refl{\alpha_2}\refl{\alpha_3}\refl{\alpha_4}\refl{\alpha_5}\refl{\alpha_6},
\\
w_6=\refl{\alpha_6}\refl{\alpha_5}\refl{\alpha_4}\refl{\alpha_2}\refl{\alpha_3}\refl{\alpha_1}\refl{\alpha_4}\refl{\alpha_3}
    \refl{\alpha_5}\refl{\alpha_4}\refl{\alpha_2}\refl{\alpha_6}\refl{\alpha_5}\refl{\alpha_4}\refl{\alpha_3}\refl{\alpha_1}.
\end{gather*}

\begin{figure}[h] \centering
\includegraphics{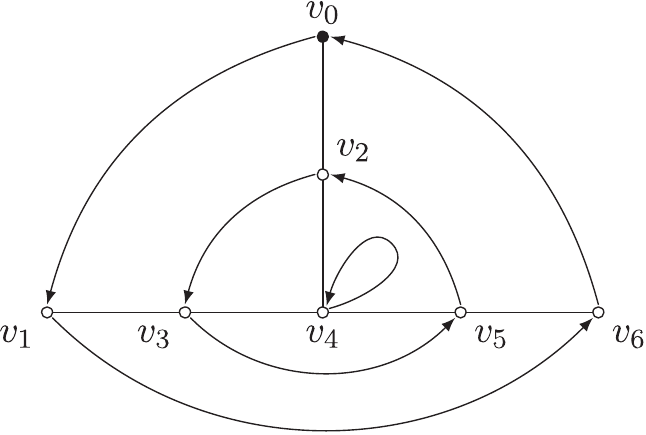}

\caption{Permutation induced by action of $\exp \lambda_1^\vee$ on the vertices of the alcove for $E_6$.}
\label{fig:e6}
\end{figure}

The permutation of the vertices induced by the action of $\exp(\lambda_1^\vee)$ (encoded by the automorphism $w_1$ of
the underlying extended Dynkin diagram) is shown schematically in Fig.~\ref{fig:e6}.

{\bf $\boldsymbol{G=E_7}$.}
Let $\mathfrak{t} \cong \mathfrak{t}^*\cong \{(x_1, \ldots, x_8) \in {\mathbb{R}}^8 \colon x_7=-x_8 \}$.
The simple roots $\alpha_1, \ldots, \alpha_7$ and highest root $\tilde\alpha$ determine the half-spaces whose
intersection is the alcove $\Delta \subset \mathfrak{t}$.
The vertices of~$\Delta$ (opposite the facets parallel to the corresponding root hyperplanes) are given in
Table~\ref{table:alcE7}.

\begin{table}[h] \centering\renewcommand{\arraystretch}{1.3}
\caption{Alcove data for $E_7$.}
\label{table:alcE7}
\vspace{1mm}

\begin{tabular}{ll}
\hline
Simple or dominant root & Opposite vertex
\\
\hline
$\alpha_1=\big(\frac{1}{2},-\frac{1}{2},-\frac{1}{2},-\frac{1}{2},-\frac{1}{2},-\frac{1}{2},-\frac{1}{2},\frac{1}{2}\big)$ &
$v_1=\big(0,0,0,0,0,0,-\frac{1}{2},\frac{1}{2}\big)$
\\
$\alpha_2=(1,1,0,0,0,0,0,0)$ &
$v_2=\big(\frac{1}{4},\frac{1}{4},\frac{1}{4},\frac{1}{4},\frac{1}{4},\frac{1}{4},-\frac{1}{2},\frac{1}{2}\big)$
\\
$\alpha_3=(-1,1,0,0,0,0,0,0)$ &
$v_3=\big({-}\frac{1}{6},\frac{1}{6},\frac{1}{6},\frac{1}{6},\frac{1}{6},\frac{1}{6},-\frac{1}{2},\frac{1}{2}\big)$
\\
$\alpha_4=(0,-1,1,0,0,0,0,0)$ & $v_4=\big(0,0,\frac{1}{4},\frac{1}{4},\frac{1}{4},\frac{1}{4},-\frac{1}{2},\frac{1}{2}\big)$
\\
$\alpha_5=(0,0,-1,1,0,0,0,0)$ & $v_5=\big(0,0,0,\frac{1}{3},\frac{1}{3},\frac{1}{3},-\frac{1}{2},\frac{1}{2}\big)$
\\
$\alpha_6=(0,0,0,-1,1,0,0,0)$ & $v_6=\big(0,0,0,0,\frac{1}{2},\frac{1}{2},-\frac{1}{2},\frac{1}{2}\big)$
\\
$\alpha_7=(0,0,0,0,-1,1,0,0)$ & $v_7=\big(0,0,0,0,0,1,-\frac{1}{2},\frac{1}{2}\big)$
\\
$\tilde{\alpha}=(0,0,0,0,0,0,-1,1)$ & $v_0=0$
\\
\hline
\end{tabular}
\end{table}

The non-zero element of the centre $Z(E_7) \cong {\mathbb{Z}}/2{\mathbb{Z}}$ is given by ($\exp$ of) the minimal
dominant coweight $\lambda_7^\vee = e_6+\frac{1}{2}(e_8-e_7)$.
The corresponding element $w_7$ of the Weyl group (as in Proposition~\ref{prop:action}), inducing an automorphism of the
extended Dynkin diagram is
\begin{gather*}
w_7=\refl{\alpha_7}\refl{\alpha_6}\refl{\alpha_5}\refl{\alpha_4}\refl{\alpha_2}\refl{\alpha_3}\refl{\alpha_1}\refl{\alpha_4}\refl{\alpha_3}
    \refl{\alpha_5}\refl{\alpha_4}\refl{\alpha_2}\refl{\alpha_6}\refl{\alpha_5}\refl{\alpha_4}\refl{\alpha_3}\refl{\alpha_1}\\
\hphantom{w_7=}{}\times
    \refl{\alpha_7}\refl{\alpha_6}\refl{\alpha_5}\refl{\alpha_4}\refl{\alpha_2}\refl{\alpha_3}\refl{\alpha_4}\refl{\alpha_5}\refl{\alpha_6}\refl{\alpha_7}.
\end{gather*}

\begin{figure}[h] \centering
\includegraphics{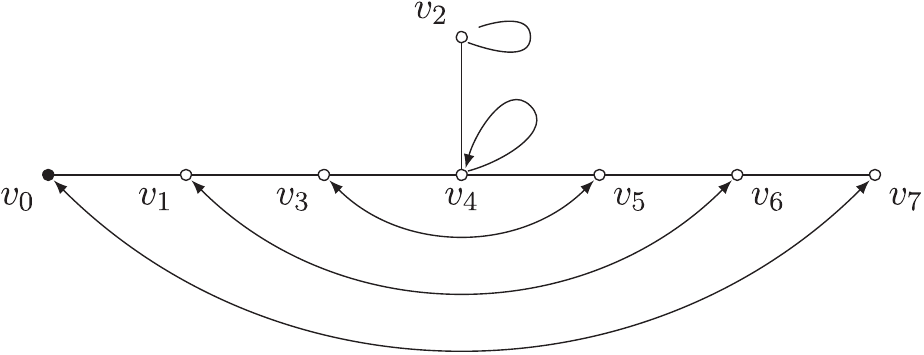}

\caption{Permutation induced by action of $\exp \lambda_7^\vee$ on the vertices of the alcove for $E_7$.}
\label{fig:e7}
\end{figure}

The permutation of the vertices induced by the action of $\exp(\lambda_7^\vee)$ (encoded by the automorphism $w_7$ of
the underlying extended Dynkin diagram) is shown schematically in Fig.~\ref{fig:e7}.

\subsection*{Acknowledgements}

A portion of this work is a~revised version of (previously unpublished) results from the author's Ph.D.\ Thesis~\cite{mythesis}, supervised by E.~Meinrenken and P.~Selick.
I~remain grateful for their guidance and support.
My sincere thanks as well to the referees for their helpful suggestions for clarif\/ication.

\pdfbookmark[1]{References}{ref}
\LastPageEnding

\end{document}